\documentclass[lettersize,journal]{IEEEtran}
\usepackage{amsmath,amsfonts}
\usepackage{algorithmic}
\usepackage{algorithm}
\usepackage{array}
\usepackage[caption=false,font=normalsize,labelfont=sf,textfont=sf]{subfig}
\usepackage{textcomp}
\usepackage{stfloats}
\usepackage{url}
\usepackage{verbatim}
\usepackage{graphicx}
\usepackage{cite}
\usepackage[switch]{lineno}
\usepackage{booktabs}
\usepackage{hyperref}
\hypersetup{hidelinks}
\usepackage{tikz,filecontents,array}
\usetikzlibrary{arrows,shapes,chains}
\usepackage{soul,color,stfloats}
\usepackage{autobreak}
\usepackage{nomencl}
\makenomenclature

\graphicspath{{./images/}}

\begin{document}

\title{A Model-Adaptive Clustering Method for Low-Carbon Energy System Optimization}

\author{
    Yuheng~Zhang,~\IEEEmembership{Student~Member,~IEEE,} Vivian~Cheng, Dharik~S.~Mallapragada, \\Jie Song,~\IEEEmembership{Senior Member,~IEEE,} Guannan~He,~\IEEEmembership{Member,~IEEE}
\thanks{Yuheng Zhang is with Department of Energy and Resources Engineering, College of Engineering, Peking University, Beijing 100871, China.}
\thanks{Vivian Cheng is with Department of Civil and Environmental Engineering, Massachusetts Institute of Technology, Cambridge, MA 02139, USA.}
\thanks{Dharik Mallapragada is with MIT Energy Initiative, Massachusetts Institute of Technology, Cambridge, MA 02139, USA.}
\thanks{Jie Song and Guannan He are with Department of Industrial Engineering and Management, College of Engineering, Peking University, Beijing 100871, China. Guannan He was also with MIT Energy Initiative, Massachusetts Institute of Technology, Cambridge, MA 02139, USA. (email: gnhe@pku.edu.cn)}
}
\markboth{}%
{Shell \MakeLowercase{\textit{et al.}}: A Sample Article Using IEEEtran.cls for IEEE Journals}


\maketitle

\begin{abstract}
Intermittent renewable energy resources like wind and solar pose great uncertainty of multiple time scales, from minutes to years, on the design and operation of power systems. Energy system optimization models have been developed to find the least-cost solution to matching the uncertainty with flexibility resources. However, input data that capture such multi-time-scale uncertainty are characterized with a long time horizon and bring great difficulty to solving the optimization model. Here we propose an adaptive clustering method based on the decision variables of optimization model to alleviate the computational complexity, in which the energy system is optimized over selected representative time periods instead of the full time horizon. The proposed clustering method is adaptive to various energy system optimization models or settings, because it extracts features from the optimization models. Results show that the proposed clustering method can significantly lower the error in approximating the solution with the full time horizon, compared to traditional clustering methods.
\end{abstract}


\nomenclature[01]{\(i, R\)}{Index and set of renewable energy resources}
\nomenclature[02]{\(s, S\)}{Index and set of storage}
\nomenclature[03]{\(j ,J\)}{Index and set of thermal generation}
\nomenclature[04]{\(t, \tau\)}{Index of time intervals}
\nomenclature[05]{\(T\)}{Set of time intervals}
\nomenclature[06]{\(c_{\text{DEG},s}\)}{Unit operation cost for storage $s$ to charge and discharge (\$/MWh)}
\nomenclature[07]{\(c_{\text{OP},j}\)}{Unit operation cost for thermal plant of type $j$ (\$/MWh)}
\nomenclature[08]{\(c_{\text{UpDn},j}\)}{Unit operation cost for each action of thermal start-up and shut-down of type $j$ (\$)}
\nomenclature[09]{\(c_{\text{ENE},s}^{\text{INV}}\)}{Investment cost for storage energy capacity of type $s$ (\$/MWh)}
\nomenclature[10]{\(c_{\text{POW},s}^{\text{INV}}\)}{Investment cost for storage power capacity of type $s$ (\$/MW)}
\nomenclature[11]{\(c_{\text{IRE},i}^{\text{INV}}\)}{Investment cost for intermittent renewable energy capacity of type $i$ (\$/MWh)}
\nomenclature[12]{\(c_{\text{THE},j}^{\text{INV}}\)}{Investment cost for thermal capacity of type $j$ (\$/MWh)}
\nomenclature[13]{\(x_{\text{DIS},s,t}\)}{Amount of discharged electricity from storage system $s$ at time $t$ (MWh)}
\nomenclature[14]{\(x_{\text{CHA},s,t}\)}{Amount of charged electricity from storage system $s$ at time $t$ (MWh)}
\nomenclature[15]{\(x_{\text{THE},j,t}\)}{Amount of electricity generated from thermal resource $j$ at time $t$ (MWh)}
\nomenclature[16]{\(y_{\text{ENE},s}\)}{Capacity of storage energy of type $s$ (MWh)}
\nomenclature[17]{\(y_{\text{POW},s}\)}{Capacity of storage power of type $s$ (MW)}
\nomenclature[18]{\(y_{\text{IRE},i}\)}{Capacity of intermittent renewable energy of type $i$ (MWh)}
\nomenclature[19]{\(y_{\text{THE},j}\)}{Capacity of thermal plant energy of type $j$ (MWh)}
\nomenclature[20]{\(n_{j,t}\)}{Number of online thermal generation units of type $j$ at time $t$}
\nomenclature[21]{\(n_{\text{UP},j,t}\)}{Number of start-up thermal generation units of type $j$ at time $t$}
\nomenclature[22]{\(n_{\text{DN},j,t}\)}{Number of shut-down thermal generation units of type $j$ at time $t$}
\nomenclature[23]{\(\text{N}_{j}\)}{Total number of thermal plant of type $j$}
\nomenclature[24]{\(\text{L}_{t}\)}{Demand at time $t$ (MWh)}
\nomenclature[25]{\(\xi_{\text{min}}\)}{Minimum output percentage of thermal generation units}
\nomenclature[26]{\(\xi_{\text{max}}\)}{Maximum output percentage of thermal generation units}
\nomenclature[27]{\(A_{i,t}\)}{Generation profile of intermittent renewable resource $i$ at time $t$}
\nomenclature[28]{\(w_{t}\)}{Amount of renewable resources curtailment at time $t$ (MWh)}
\nomenclature[29]{\(E_{s,t}\)}{Energy storage level of storage type $s$ at time $t$ (MWh)}
\nomenclature[30]{\(\eta\)}{Efficiency of storage charging and discharging}
\nomenclature[31]{\(\text{R}\)}{Total percentage of renewable energy in demand}
\printnomenclature
\section{Introduction}
    \IEEEPARstart{I}{n} order to keep the 1.5°C target of the Paris Agreement in reach, human beings have to act immediately and more concretely on decarbonization \cite{allen2019technical}. Countries over the world have set many targets of energy transition to lower carbon emissions. For example, in the Energy Roadmap 2050, the European Union (EU) commits itself to reducing greenhouse gas (GHG) emissions to 80-95\% below 1990 levels and realizing carbon neutrality by 2050 \cite{RN434}. China has announced its aims to reach peak carbon emissions before 2030 and carbon neutrality by 2060. To fulfill the targets, the share of fossil sources in China's energy sector needs to be reduced to less than 20\% \cite{davidson2021policies}. The USA has committed itself to the target of reducing carbon emissions to half by 2030, compared to 2005 levels, and achieving net zero emissions no later than 2050. Before the official announcement came out \cite{sheet2021president}, California state had focused on direct air carbon capture \cite{marcucci2017road} and intended to reach carbon neutrality by 2045 \cite{wheeler2017carbon}.
    
    \IEEEpubidadjcol
    
\begin{figure*}[!tbh]
    \centering
    \includegraphics[width=\textwidth]{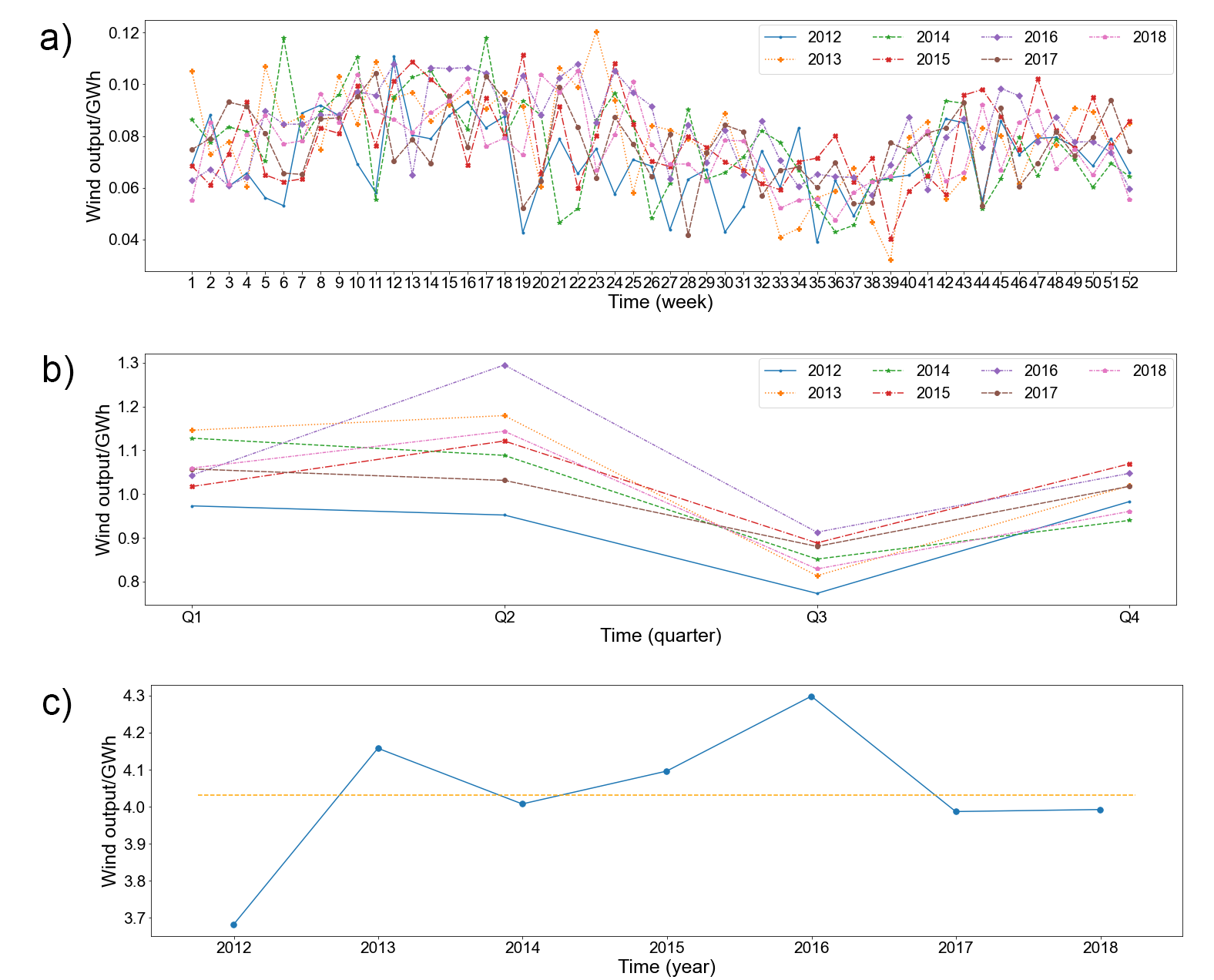}
    \caption{Multi-time-scale variability of renewable energy generation: a) weekly, b) seasonal, and c) yearly}
    \label{IRE_var}
\end{figure*}

    Intermittent renewable energy (IRE) resources are getting more and more deployed worldwide due to their sustainability and potential to meet energy demands with zero or near zero emissions of both air pollutants and GHGs \cite{RN366}. According to the International Energy Agency, IRE resources in 2020 supplied 28\% of the total world energy demand \cite{iea2020global}, reaching a 45\% increase in global capacity growth, and this share is expected to increase very significantly (30–80\% in 2100) \cite{fridleifsson2001geothermal}. However, the energy generation from IRE is variable and uncertain, in contrast to conventional thermal generation \cite{bessa2014handling}. There is not only hourly variability, but also weekly, seasonal, and yearly variability in IRE generation, as shown in Figure \ref{IRE_var}. With such variability and uncertainty, when the IRE rapidly penetrates into power grids, many grid integration problems could arise including curtailment \cite{golden2015curtailment} and over-generation \cite{denholm2015overgeneration}.
    
    To cope with the challenges, greater investments in flexibility resources like energy storage technologies are necessary. Energy storage can provide multiple values to energy systems, such as accommodating increased penetrations of IRE \cite{RN428}, load leveling and peak shaving \cite{oudalov2007sizing}, frequency regulation \cite{stroe2016operation}, damping energy oscillations \cite{zhu2018optimization}, and improving power quality and reliability \cite{ribeiro2001energy}. The global energy storage capacity is expected to triple to 181-421 GWh by 2030 \cite{ralon2017electricity}. To optimally match flexibility to variability in energy systems, energy system optimization models (ESOMs) are commonly used, which determine the least-cost sizes and scheduling of generation resources such as IRE, energy storage, and thermal power plants. When involved with integer variables to denote the status of thermal generation units or the decision of transmission line expansion, these models are in the form of mixed integer linear programming (MILP) \cite{alizadeh2011reliability,haghighat2018bilevel,khodaei2012coordination,castro2018expanding}. In order to capture the multi-time-scale balance of variability and flexibility in future energy systems, incorporating long-term data of operation from IRE and storage systems is necessary, especially if multi-stage decarbonization policies/pathways towards the mid of the century need to be modeled. Such models are also called capacity expansion \cite{RN429} problems, whose focus is on how to optimally invest new power and storage capacities for future energy supply requirements.

    For long-term ESOMs in the form of MILP, the computational complexity increases dramatically as the time horizon expands. Time aggregation \cite{RN185} or temporal clustering \cite{RN192} techniques have been applied to reduce the complexity of these models in such situations. The main concept is to represent the full time horizon with some representative time slices, selected through finding patterns or clusters of similar supply and demand in the full time series. Many time slicing/aggregation methods have been developed. Following Haydt et al. \cite{haydt2011relevance}, two commonly used methods are referred to as the integral method and the semi-dynamic method. In the integral method, typically 5-10 time slices are used to distinguish between different load levels across the time horizon and each time slice represents an average load level during a certain fraction of the full horizon \cite{haydt2011relevance}. The semi-dynamic method is the most used time slicing method \cite{pina2013high,welsch2015supporting}, which disaggregates time horizon into different seasons, days of the week and diurnal periods under the assumption that the time series data change on the basis of seasonal, weekly and daily frequency. More advanced methods are in the clustering family. Nahmmacher et al. \cite{nahmmacher2016carpe} used a hierarchical algorithm with Ward linkage to perform clustering on the LIMES-EU model and showed that a small number of representative days developed in this way are sufficient to reflect the characteristic fluctuations of the input data. Zatti et al. \cite{RN330} proposed a MILP clustering model, named k-MILP, devised to find at the same time the most representative days of the year and the extreme days. Vitali et al. \cite{RN420} proposed a clustering method considering temporal correlation as an important feature in time aggregation to offer better representation in capacity expansion problems including storage and renewable energy sources, which keeps chronology in time series. All of the existing time aggregation methods examine similarity from the time series input data themselves. However, it is neglected in the studies that whether similar data profiles have similar impacts on the ESOM results is uncertain, as the ESOM is usually a non-linear mapping from inputs to solutions. Therefore, traditional time aggregation methods based on similarity of raw data only may not guarantee high approximation performance for long-term ESOMs. 
    
    To fill the gap, here we propose a model-adaptive clustering-based time aggregation method that combines model characteristics and input data variations in representative time period selection. Features used in clustering are extracted from processing the raw input data with an ESOM and contain similarity information from the perspective of the ESOM. For different ESOMs or policy settings, this method can adaptively select representative time periods, which can be then used to approximate the full time horizon and reduce computational complexity with lower approximation errors than traditional methods. The main contributions are three-fold: 
    1) we propose to extract features through dimension reduction using the ESOM for representative time period selection;  
    2) we find the best clustering configurations for the proposed model-adaptive time aggregation method;
    3) we show how the proposed method is adaptive to different policy settings, which explains why the approximation errors of the proposed method are lower than traditional methods.
    
    The remainder of this article is organized as follows: Section 2 describes and establishes a case ESOM including thermal generation and renewable energy resources as well as storage systems. Section 3 explains in detail how to implement our proposed adaptive clustering method. Numeric results are presented in Section 4. Section 5 draws conclusions.

\section{Energy System Optimization Model}
\label{Energy System Optimization Model}
    Bottom-up, long-term energy system planning models or capacity expansion models (CEMs) are frequently used to analyze pathways for the transition of the energy/electrical power system and to inform policy design \cite{poncelet2020unit}. Without loss of generality, we design a case ESOM that contains different technologies including thermal, wind, and solar generators as well as energy storage systems, to compare the performance of different time aggregation methods and validate our proposed method. 

\subsection{Objective function}
     The objective function consists of variable costs in operation and fixed costs in investment denoted by $C_{\text{VAR}}$ and $C_{\text{FIX}}$, as in equation (\ref{OB}). The decision variables are the capacities and operational schedules of each technology. For many storage technologies, the power (conversion capability) and energy (storage volume) capacities can be designed independently, denoted by $y_{\text{POW},s}$ and $y_{\text{ENE},s}$, respectively. 

\begin{equation}
    \label{OB}
    OB = C_{\text{VAR}} + C_{\text{FIX}}
\end{equation}
\begin{equation}
        C_{\text{VAR}} = C_{\text{DEG}} + C_{\text{OP}} + C_{\text{UpDn}}
    \label{CVAR}
\end{equation}
\begin{equation}
    C_{\text{DEG}} = \sum_{t \in \mathbb{T}} \sum_{s \in \mathbb{S}} (x_{\text{DIS},s,t} + x_{\text{CHA},s,t}) c_{\text{DEG},s}
    \label{CDEG}
\end{equation}
\begin{equation}
    C_{\text{OP}} = \sum_{t \in \mathbb{T}} \sum_{j \in \mathbb{J}} x_{\text{THE},j,t} c_{\text{OP},j}
    \label{COP}
\end{equation}
\begin{equation}
    C_{\text{UpDn}} = \sum_{t \in \mathbb{T}} \sum_{j \in \mathbb{J}} (n_{\text{UP},j,t} + n_{\text{DN},j,t}) c_{\text{UpDn},j}
    \label{CUPDOWN}
\end{equation}
\begin{equation}
    \label{CFIX}
    \begin{split}
        C_{\text{FIX}} = & \sum_{s \in \mathbb{S}} (y_{\text{ENE},s} c_{\text{ENE},s}^{\text{INV}} + y_{\text{POW},s} c_{\text{POW},s}^{\text{INV}}) + \\ & \sum_{i \in \mathbb{R}} y_{\text{IRE},i} c_{\text{IRE},i}^{\text{INV}} + \\ & \sum_{j \in \mathbb{J}} y_{\text{THE},j} c_{\text{THE},j}^{\text{INV}}\text{N}_{j}
    \end{split}
\end{equation}

    Variable costs are mainly from discharge and charge actions in the storage systems, denoted by $C_{\text{DEG}}$, and thermal operation, which consists of two parts: operation cost and startup and shutdown costs, which are denoted by $C_{\text{OP}}$ and $C_{\text{UpDn}}$, respectively. Fixed costs consist of investment costs from different technologies considered in system including energy capacity, power capacity for storage system, multi renewable capacity and thermal capacity, as in equation (\ref{CVAR}).
    
\subsection{Energy balance}

    For each IRE resource considered in the model, the capacity $y_{\text{IRE},i}$ indicates its maximum possible generation, and its profile is represented by $A_{i,t}$. When an IRE resource is over-abundant, curtailment might be needed and represented by $w_{t}$. Considering all energy inflow and outflow, the energy balance can be modelled as:  

\begin{equation}
     \sum_{i \in \mathbb{R}} y_{\text{IRE},i} A_{i,t} -w_t + x_{\text{THE},j,t}  + x_{\text{DIS},s,t}-x_{\text{CHA},s,t} = \text{L}_t
    \label{balance}
\end{equation}

    In equation (\ref{balance}), $y_{\text{IRE},i} A_{i,t}$ represents the output of IRE resource $i$ at time $t$ in the form of capacity multiplied by its generation profile; $x_{\text{THE},j,t}$ denotes the output of thermal generator $j$ at time $t$; $x_{\text{DIS},s,t}$ and $x_{\text{CHA},s,t}$ stand for discharged and charged energy from storage systems at time $t$, respectively; $w_t$ denotes the amount of IRE curtailment; and $L_t$ denotes power demand at time $t$. Each term of equation (\ref{balance}) should be positive.

\subsection{Thermal generation}

    The minimum and maximum generation requirements of thermal power plants when they are online are modelled as:

\begin{equation}
    x_{\text{THE},j,t}\geq \xi_{min} y_{\text{THE},j,t} n_{j,t}
\end{equation}
\begin{equation}
    x_{\text{THE},j,t}\leq \xi_{max} y_{\text{THE},j,t} n_{j,t}
\end{equation}
\begin{equation}
    \xi \in [0,1], \xi_{min} \leq \xi_{max}
\end{equation}
where $\xi_{min}$ and $\xi_{max}$ denote the minimum and maximum output, respectively, as percentages of the full capacity $y_{\text{THE},j,t}$. $n_{j,t}$ denotes the number of online units. The unit commitment of thermal generators are modelled as:

\begin{equation}
    \label{start_shut}
    n_{j,t} - n_{j,t-1} = n_{\text{UP},j,t} - n_{\text{DN},j,t}
\end{equation}

\begin{equation}
    0 \leq n_{j,t}, n_{\text{UP},j,t}, n_{\text{DN},j,t} \leq N_{j}
\end{equation}
\begin{equation}
    n_{j,t}, n_{\text{UP},j,t}, n_{\text{DN},j,t},y_{\text{THE},j} \in \mathbb{N}
\end{equation}

    Namely the subscripts 'UP' and 'DN' denote the start-up and shut-down actions of thermal generation units, respectively. The minimum up/down time constraints respectively require units to remain online/offline for a minimum period of time after starting up/shutting down. These constraints are formulated as equations (\ref{start_shut})-(\ref{down_time}), where the minimum up and down time are denoted by $\tau^{\text{UP}}$ and $\tau^{\text{DN}}$.

\begin{equation}
    \label{up_time}
    n_{j,t} \geq \sum_{\tau=t-\tau^\text{UP}}^t n_{\text{UP},j,t}
\end{equation}
\begin{equation}
    \label{down_time}
    N_{j} - n_{j,t} \geq \sum_{\tau=t-\tau^\text{DN}}^t n_{\text{DN},j,t}
\end{equation}
\subsection{Storage system}

    The state of charge and operational limit constraints for energy storage systems are modelled as:

\begin{equation}
    E_{s,t} - E_{s,t-1} =  x_{\text{CHA},s,t}\eta - \frac{x_{\text{DIS},s,t}}{\eta}
\end{equation}
\begin{equation}
    0 \leq E_{s,t} \leq y_{\text{ENE},s}
\end{equation}
\begin{equation}
    0 \leq x_{\text{DIS},s,t}, x_{\text{CHA},s,t} \leq y_{\text{POW},s}
\end{equation}
\begin{equation}
    E_{s,t}, x_{\text{DIS},s,t}, x_{\text{CHA},s,t} \in \mathbb{R}^{+}
\end{equation}

    $E_{s,t}$ is the energy level at time slot $t$ in the storage system for each type of storage device. $\eta$ is efficiency of battery charging and discharging. $x_{\text{CHA},s,t}$ and $x_{\text{DIS},s,t}$ are the amount of electricity charged and discharged from storage device $s$. Storage energy level can't exceed its capacity and charging power can't overshoot power capacity.

    We model renewable portfolio standards, which are very common energy policies over the world, setting lower bounds for energy supplied by renewable energy and equivalently upper bounds for fossil fuel based generation as:

\begin{equation}
    \frac{\sum_{t=1}^{T} x_{\text{THE},j,t}}{\sum_{t=1}^{T} \text{L}_t} \leq 1 - \text{R}
\end{equation}

    This $\text{R}$ parameter controls the portfolio of dispatchable thermal generation and IRE along with the storage system. A bigger $\text{R}$ requires a larger share of IRE and storage systems.

\section{Model-Adaptive Clustering Method}
\label{Model-Adaptive Clustering Method}
    As introduced before, there is multi-time-scale variability in IRE generation, from minutes to years, and thus, solving ESOMs across a long time horizon is needed to design an energy system with high penetrations of IRE.
    
\begin{figure*}[!htb]
    \centering
    \includegraphics[width=\textwidth]{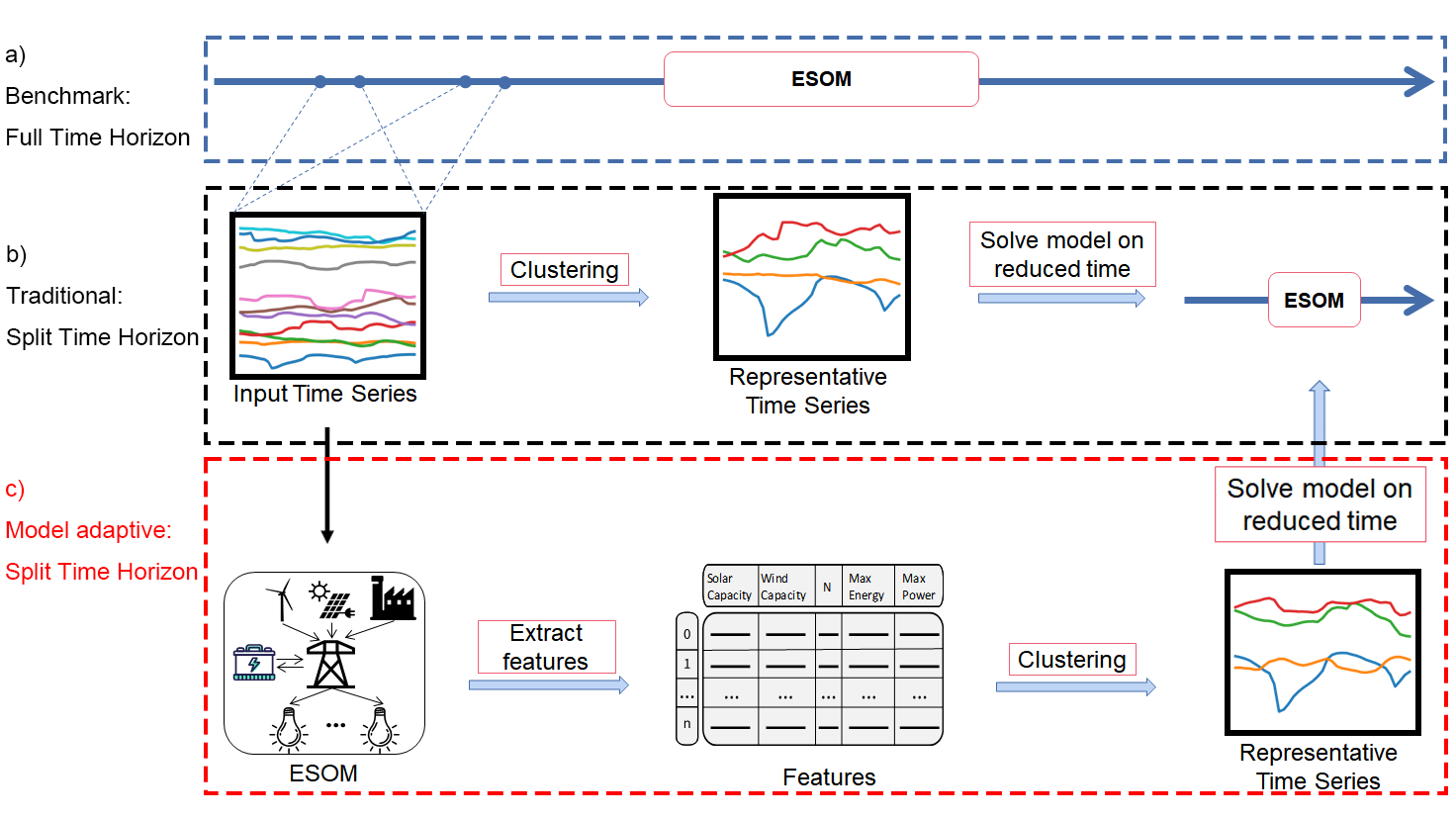}
    \caption{Schematics of time aggregation: a) benchmark: solving full time horizon model; b) traditional time aggregation: clustering on raw data; c) model adaptive time aggregation: clustering on features extracted from models}
    \label{schematic}
\end{figure*}

    In long-term ESOMs, input data including IRE sources and demand profiles span across a long time horizon. Such ESOMs have extremely high dimensions, and computational intractability is a significant issue when directly solving the models \cite{RN293}, as in Figure (\ref{schematic}).a). Time aggregation is used to address the intractability through approximation.
    
    In traditional time aggregation, as in Figure (\ref{schematic}).b), the full time horizon is first divided into smaller time slices. Then, clustering algorithms are applied to the time slices to group similar slices into clusters based on pre-defined distances between time slices, and one slice is selected to represent all slices in a cluster. Solving the ESOM over the selected slices could save computational time significantly.

    Unlike traditional time aggregation in which raw data of IRE and demand profiles are used as clustering features, our proposed model-adaptive method extracts features from ESOMs for clustering, and thus, the clustering results contain model information, as shown in Figure (\ref{schematic}).c). The detailed procedures are as follows:

    1) Choose an appropriate period length, which usually is a day or week, and split the full planning horizon $T$ into a set of such smaller periods denoted by ${T_1^{'}, T_2^{'}, \cdots, T_m^{'}}$;
    
    2) Run the ESOM of interest on smaller period $T_i^{'}, i \in 1, 2, \cdots, m$ and extract features from the model;
    
    3) On the basis of procedure 2), a mapping between period and features is established. Perform an appropriate clustering algorithm on the features and obtain representative periods along with corresponding weights. The weights denote how many periods in the whole time horizon the selected periods represent;
    
    4) Solve the ESOM on the representative periods selected.

    In adaptive time aggregation, by implementing procedure 2), the features containing both information from the ESOM and original input data variation along time are extracted for each time slice. In other words, the features are obtained through dimension reduction, using the ESOM as a non-linear transformer. These intrinsic features can be extracted from anywhere in the ESOM such as decision variables, dual variables, etc., as long as they contribute to finding the best representative time slices to approximate the benchmark solution of the full time horizon. For CEMs, specifically, the features could be planned capacities or costs of each type of generation technology. From the ESOM defined before, five decision variables are selected: solar capacity, wind capacity, storage energy capacity, storage power capacity, and thermal generator capacity. After solving the ESOM for each time slice, the features of each time slice are computed, and time slices with similar capacity planning outcomes are clustered together and represented by one time slice.

\begin{figure}
    \centering
    \includegraphics[width=3.45in]{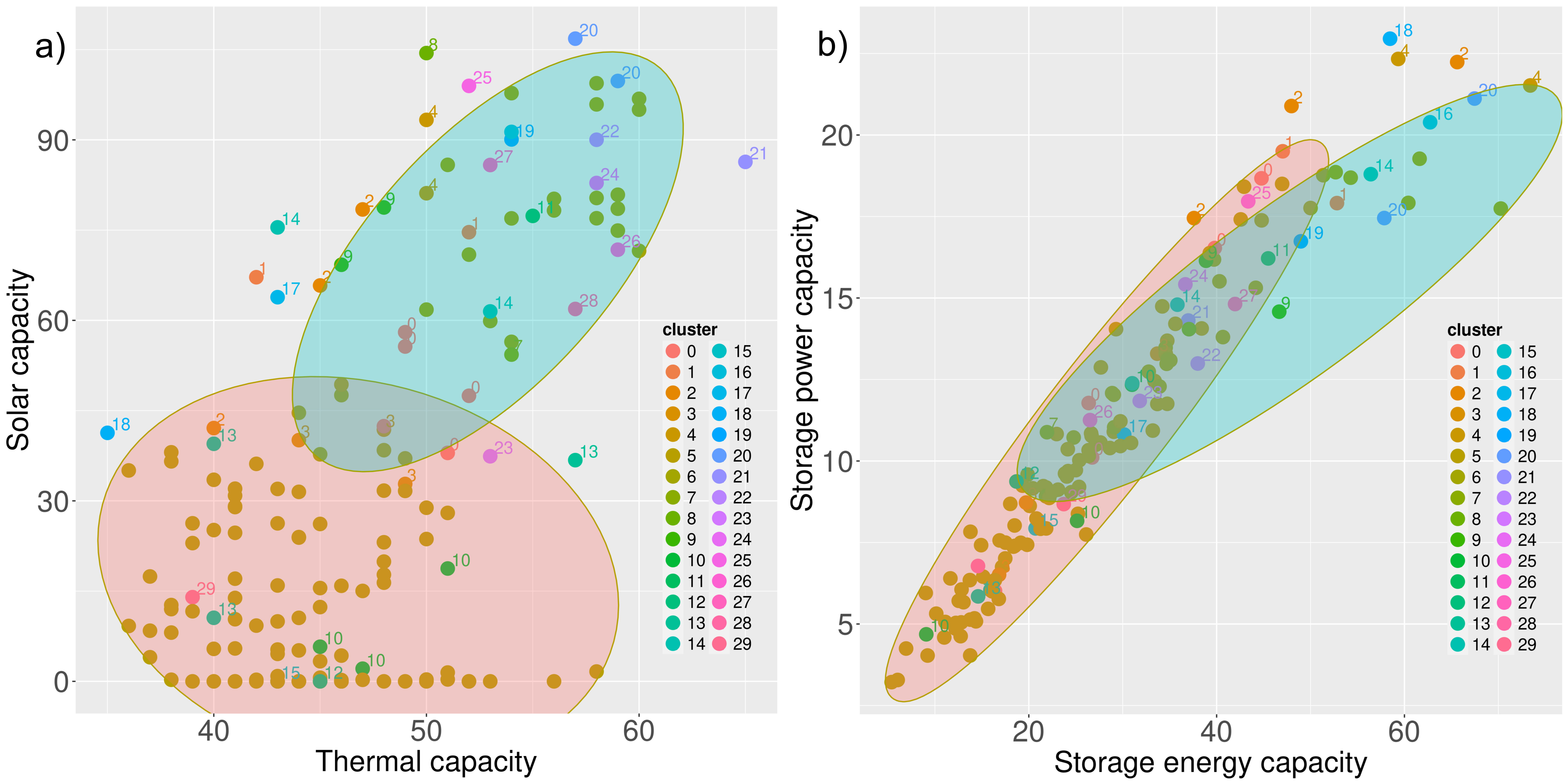}
    \caption{Features of decision variables clustering: a) thermal \& solar capacity b) storage energy and power capacity}
    \label{features}
\end{figure}

    To evaluate the performance of our proposed and traditional time aggregation methods in approximating the benchmark that solves the ESOM with the full time horizon \cite{RN411}, we use the metric of mean absolute percentage error (MAPE) in planning outcomes, defined as: 

\begin{equation}
\text{MAPE}=\frac{1}{N}\sum_{i}^{N} \frac{\left|\hat{y}_{i}-y_{i}\right|}{y_{i}}
\end{equation}
where $y_{i}$ is the $i$th capacity planning outcome of the benchmark solution, and $\hat{y}_{i}$ is the $i$th capacity planning outcome estimated using a time aggregation method, either traditional or adaptive.

\section{Case Study}
\label{Case Study}
\subsection{Case setting}

    Two case studies are conducted to validate the method, one with a linear ESOM and the other with an integer one. The integer case is a full version of the ESOM presented in Section 2 in which only one type of thermal generation is used. In the linear case, thermal generation is removed from the portfolio, and the power demand is met by wind, solar, and energy storage systems only, which models a 100\% renewable energy system. The integer variables for modelling the unit commitment of thermal generation units are therefore not needed. For the integer cases, the time horizon is reduced from 7 to 3 years, to keep the benchmark model solvable in a reasonable time frame. The code for both cases can be found in \cite{Zhang_Adaptive-clustering-for-ESOM_2021}.
    
    To capture enough variability in a single representative period, here a week's data is bundled together as a time slice, which is 168 points in the resolution of an hour. Many articles \cite{RN509, RN195} select representative days using clustering methods to ease computational burden with the planning horizon usually being no more than a year, which is far from requirements for data in long-term ESOMs, especially with renewable energy and storage systems.
    
    K-means and hierarchical clustering methods with different linkages are applied to compare how different clustering methods affect the accuracy of time aggregation. 
    In the base case, we use agglomerative clustering with single linkage to select 30 representative time slices, and this clustering setting will be justified later.  
    
    In order to generate different ESOM settings, the renewable portfolio standard is varied from 50\% to 95\%, yielding 10 different policy constrained scenarios. Other model parameters related to investment and operation are listed in Table \ref{parameters}.

\begin{table}[!htbp]
\caption{Main parameters in the ESOM.}
\centering
\resizebox{\columnwidth}{!}{
    \begin{tabular}{lll}
    \toprule
    Model Parameters &Type &Value\\
    \midrule
    Battery Energy &CAPEX &$200$ \$/kWh\\

    Battery Power &CAPEX &$70$ \$/kW\\

    Photovoltaic Panel &CAPEX &$1000$ \$/kW\\

    Wind Turbine &CAPEX &$1500$ \$/kW\\

    Thermal Plant &CAPEX &$1000$ \$/kW\\

    Battery Charing and Discharging &OPEX &50 \$/MWh\\

    Thermal Generation &OPEX &30 \$/MWh\\

    Up Time &Other &6 hrs\\

    Down Time &Other &6 hrs\\
    \bottomrule
    \end{tabular}
}

\label{parameters}
\end{table}
\subsection{Performance of adaptive clustering}

\begin{figure}[!tbh]
    \centering
    \includegraphics[scale=0.35]{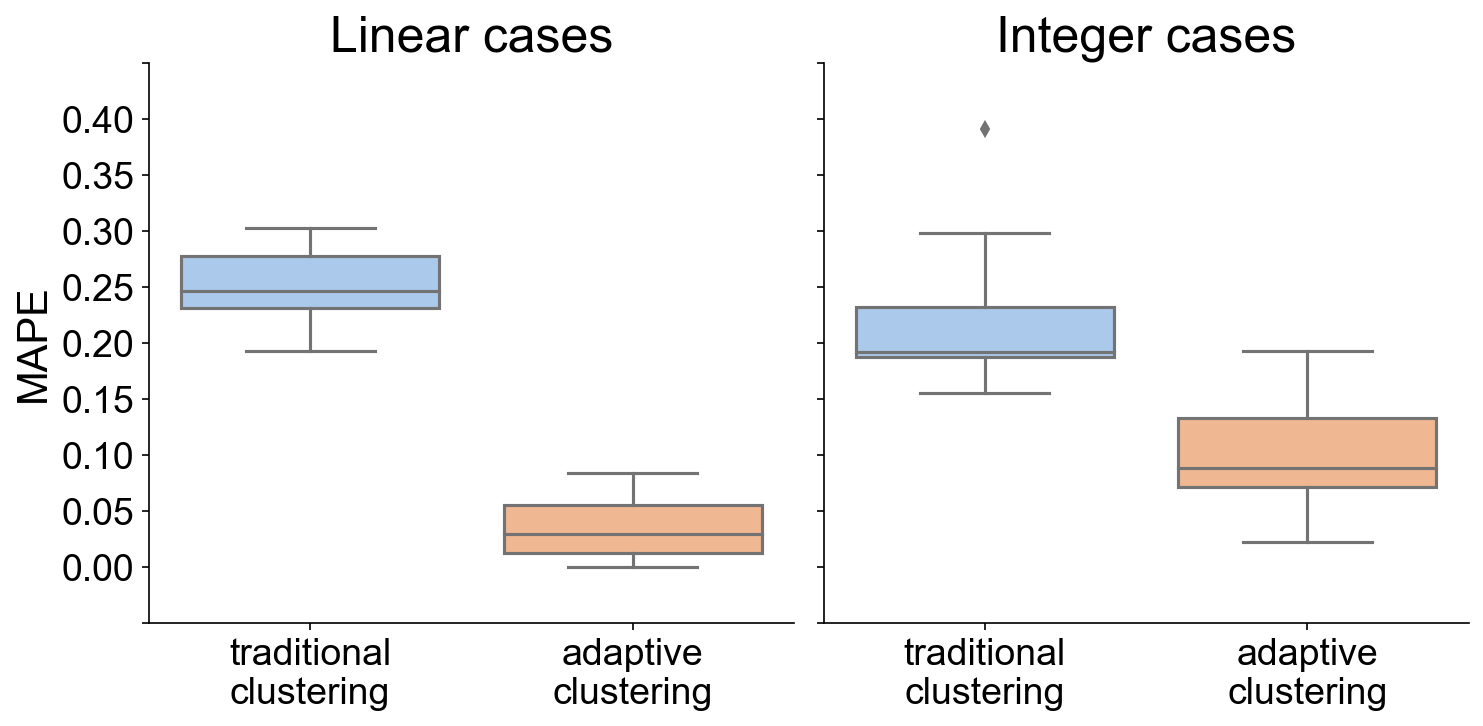}
    \caption{MAPE of decision variables of traditional and adaptive clustering with linear and integer cases. In linear cases, costs of wind and solar are altered and 16 different combinations are used to estimate average performance, while in integer cases, 10 scenarios are used where different penetration of renewable energy varying from 0.5 to 0.95 is set. Error bars indicate the variation from parameters in the model.}
    \label{mape_linear_integer}
\end{figure}

\begin{figure}[tbh]
    \centering
    \includegraphics[scale=0.25]{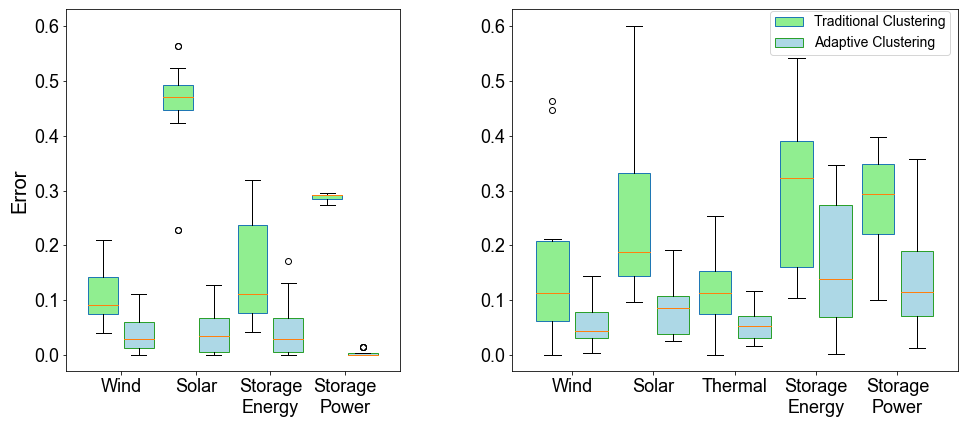}
    \caption{Errors of decision variables of traditional and adaptive clustering with linear and integer cases. In linear cases, capacity of wind, solar, storage energy and storage power are decomposed to show errors of each term, while in integer cases, thermal capacity is further modeled and measured.}
    \label{mape_linear_integer_decomposed}
\end{figure}

    In Figure \ref{mape_linear_integer}, the results of the linear cases show a 21.5\% reduction in MAPE using our proposed adaptive clustering method compared to the traditional method. The MAPE of the adaptive clustering method can be as low as 3.5\%. In the integer cases, the adaptive clustering also brings a 12.5\% decrease in MAPE. In Figure \ref{mape_linear_integer_decomposed}, the adaptive clustering method has lower errors for each of the capacity decision variables. It performs better in capturing IRE's characteristics given that in both linear and integer cases the errors on capacity of wind and solar are on average lower than half of those in traditional clustering. 
    
    In terms of solving efficiency, linear and integer cases take 335 seconds and 2 hours in solving benchmark models, respectively. While using model-adaptive clustering, it takes about 244 seconds to train and obtain representative weeks for linear cases and 245 seconds for integer ones. Considering that only three years of data are used in integer cases, it might be scaled to 572 seconds to train a seven-year case of integer model. And the elapsed time is shortened to 10 seconds in linear cases and 89 seconds in integer ones. Our proposed method could remarkably reduce the total computational time to less than five minutes, eliminating much computational burden, especially for MILP.
    
\begin{figure}[tbh]
    \centering
    \includegraphics[scale=0.36]{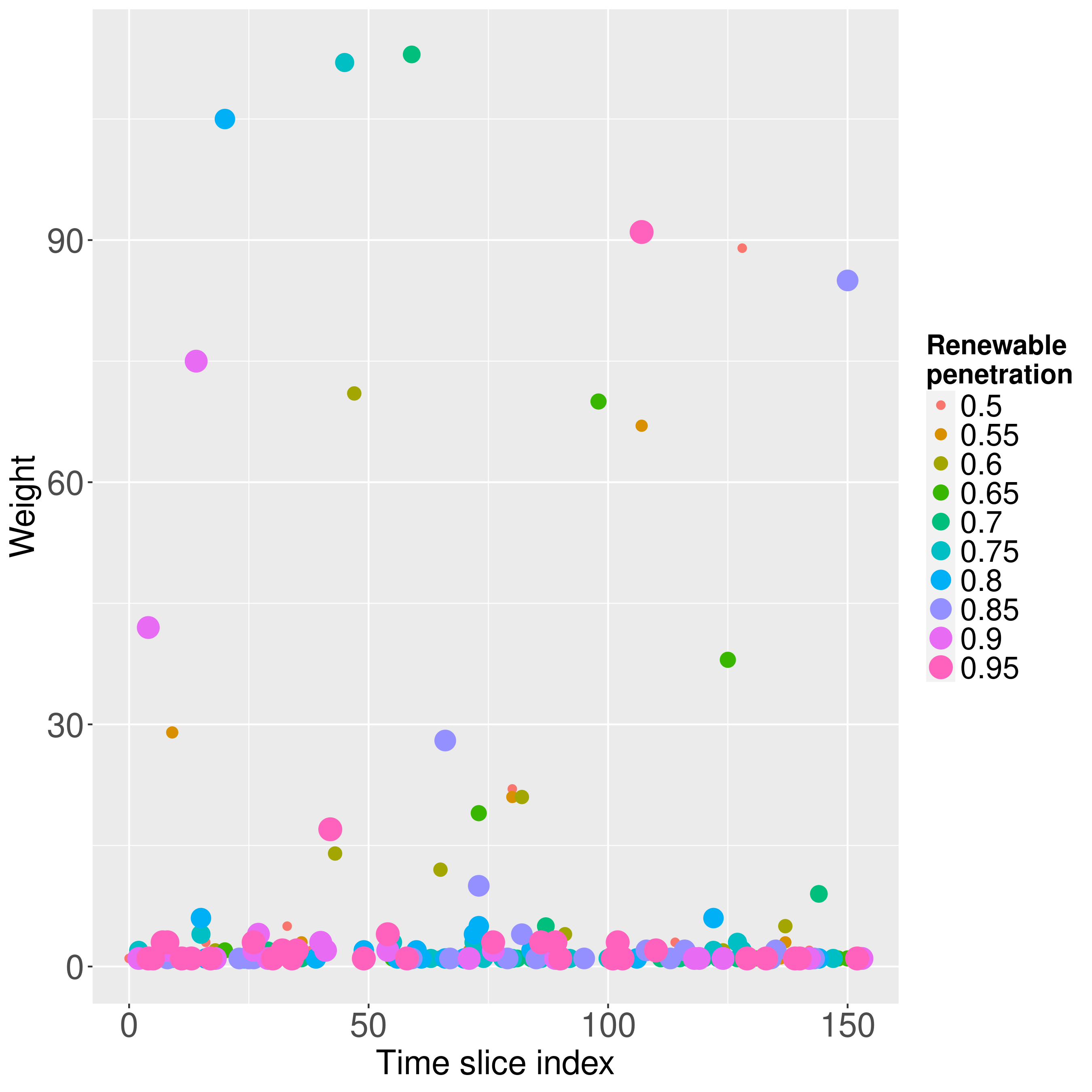}
    \caption{Selected time periods and weights with varying renewable portfolio standards.}
    \label{indices_and_weights}
\end{figure}

    Why is our proposed adaptive method effective? Intuitively, our adaptive clustering process embeds similarity information from the ESOM’s perspective, and thus the selected time periods have greater power in representing the whole time horizon. Figure \ref{indices_and_weights} shows how the representative period selection is adaptive to different model policy settings. If the weight is 0, then this time slice is not selected and is represented by another slice. If the weight is positive, this time slice is selected and represents some slices after clustering. Different colors represent the weights with different renewable portfolio policies, from 50\% to 95\%. The system setting changes when the renewable penetration changes. As shown in Figure \ref{indices_and_weights}, when the system changes, the clustering results change as well. For example, with 85\% renewable penetration, the week with the highest weight is the week number 150, while with 80\% renewable penetration, the week with the highest weight is the week number 20. This indicates that our clustering method is adaptive to the change in ESOM.

\subsection{Appropriate number of representative time periods}

\begin{figure}[H]
    \centering
    \includegraphics[width=3.45in]{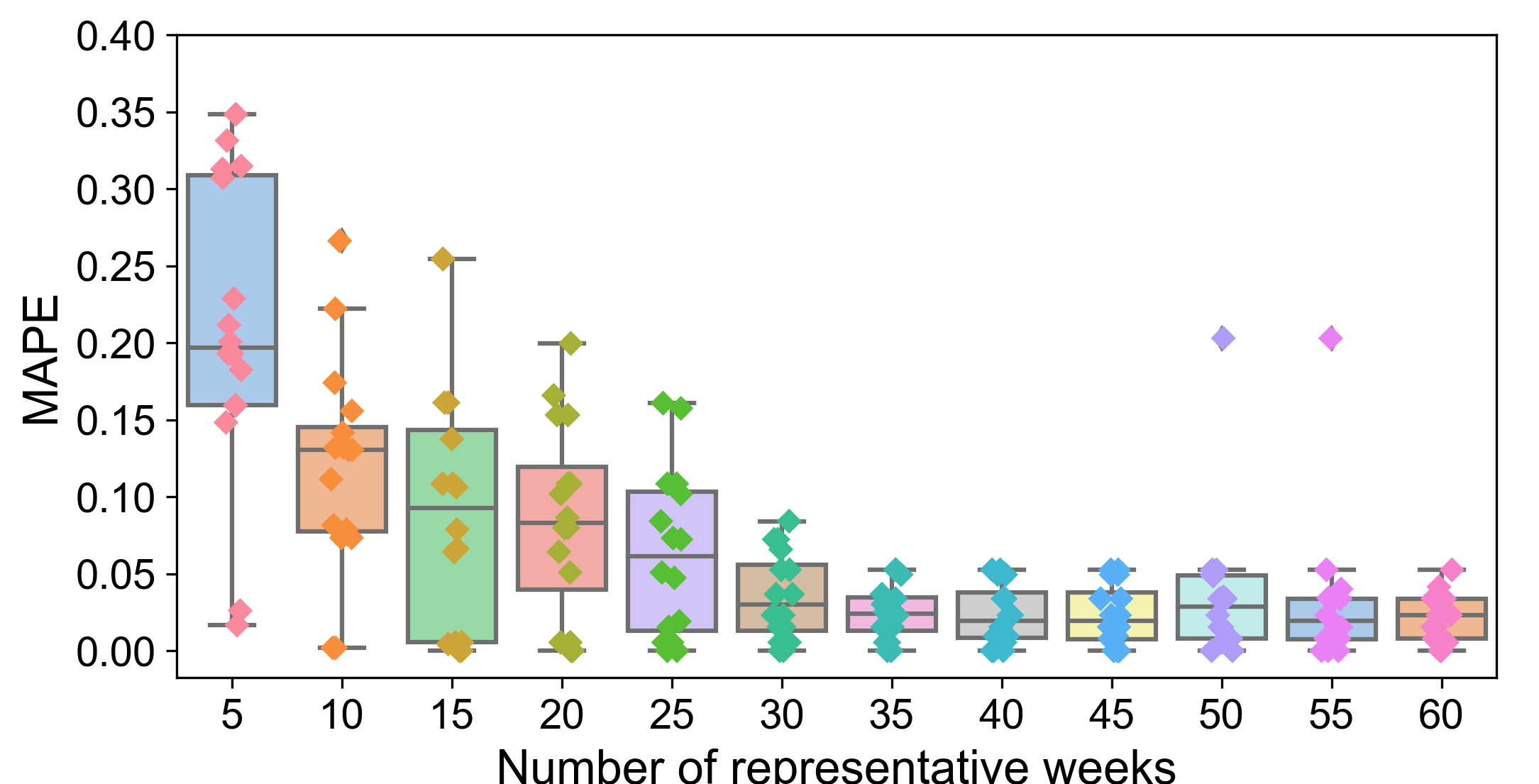}
    \caption{MAPE with different clustering numbers}
    \label{clustering_numer}
\end{figure}

    Generally we regard that to some extent, there exists a certain number of representative time periods that include enough information to approximate the full time horizon with limited errors. So with the number of clusters getting close to this level, the mean absolute error should drop to an appropriate level and the results will be more acceptable. In \cite{RN419}, Elbow’s rule \cite{RN424} is adopted to select an appropriate clustering number for picking representative periods to solve a transmission expansion planning, and a lower operating cost error below 3\% is obtained. There are many other criteria to select an appropriate clustering number \cite{RN424}. Figure \ref{clustering_numer} shows the approximation errors of adaptive time aggregation with different numbers of clusters. When the number of clusters reaches 30, the MAPE decreases slowly as the number of clusters increases. Therefore, 30 weeks seem to be enough to represent the full 7-year horizon in the case study.

    Different clustering methods perform diversely in MAPE, among which agglomerative clustering with single linkage has an average error of 3.5\%. Other methods like K-means, which is often used in representative time period selection, and other linkage modes yield errors with a mean between 10\% and 15\%. Although agglomerative clustering with average linkage sometimes can yield a lower error, its mean performance is worse than that with single linkage. These results justify our clustering setting for the base case.

\begin{figure}[H]
    \centering
    \includegraphics[width=3.45in]{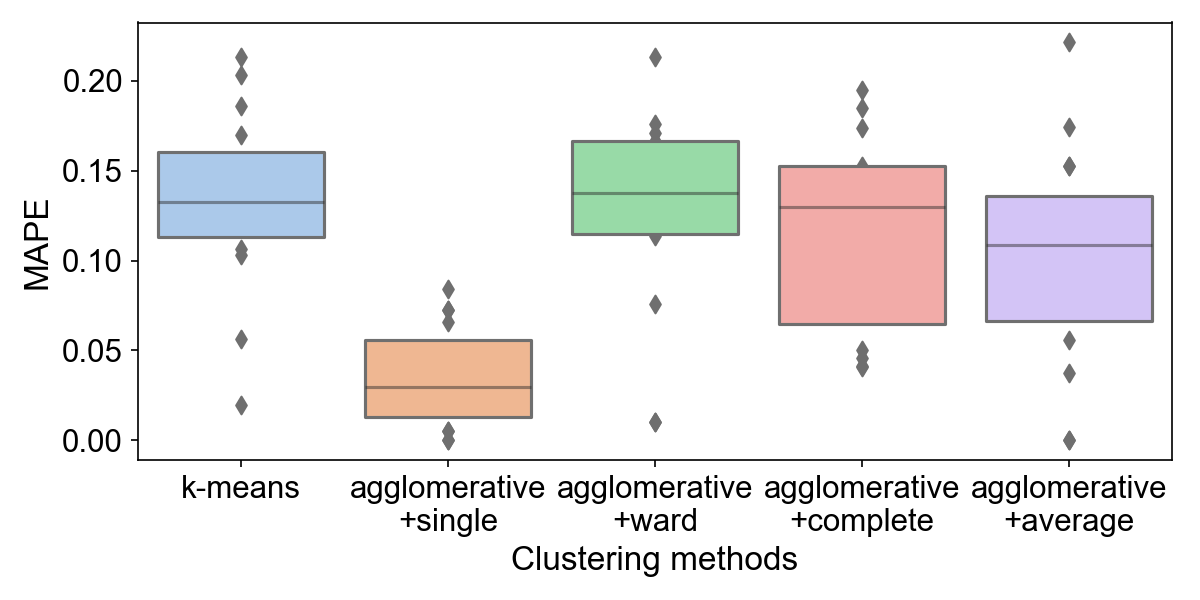}
    \caption{MAPE with different clustering methods}
    \label{clustering_method}
\end{figure}

\subsection{The smoothing effect of cluster centroid}

\begin{figure}
    \centering
    \includegraphics[scale=0.56]{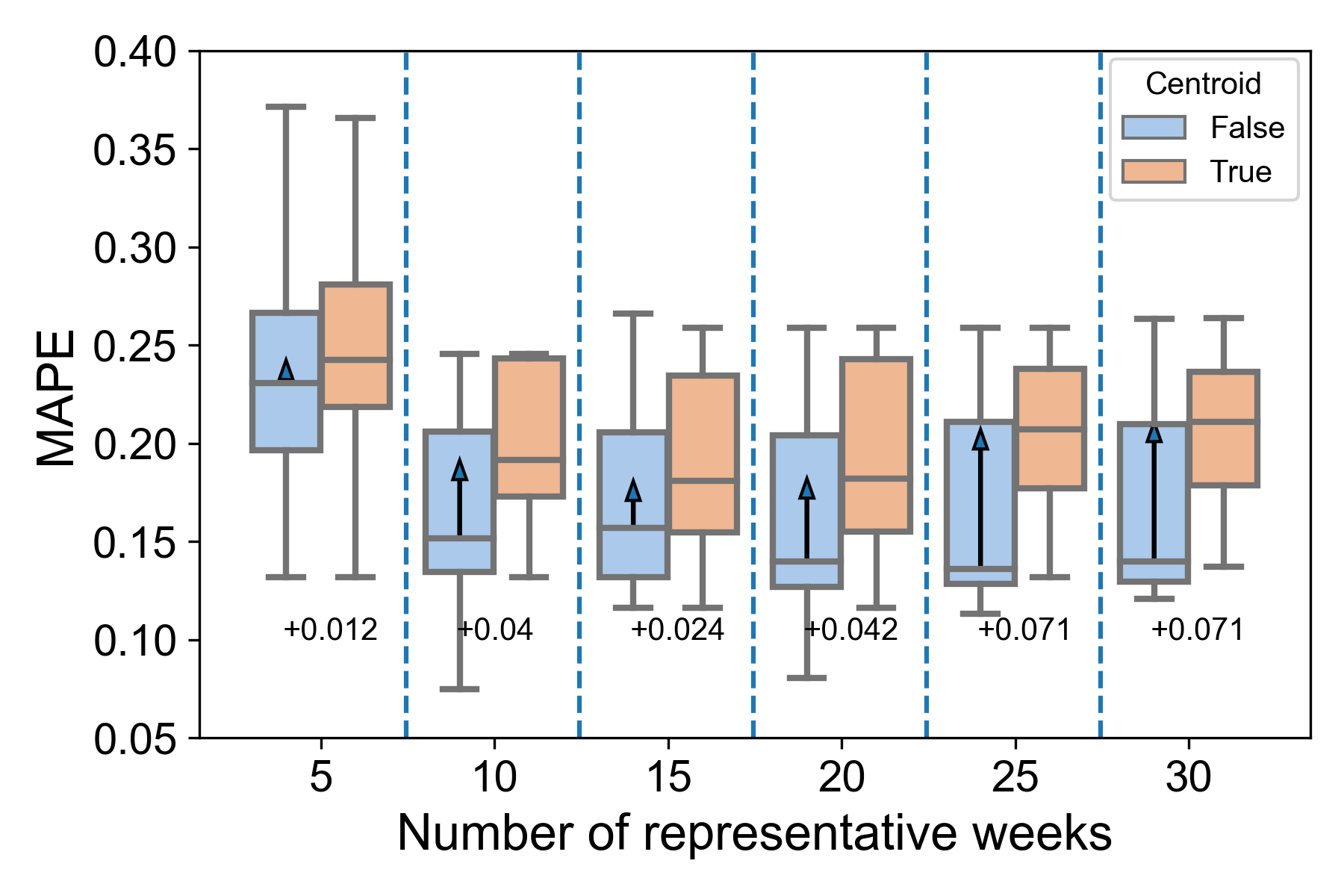}
    \caption{MAPE with different centroid setting}
    \label{centroid}
\end{figure}

\begin{figure}
    \centering
    \includegraphics[scale=0.29]{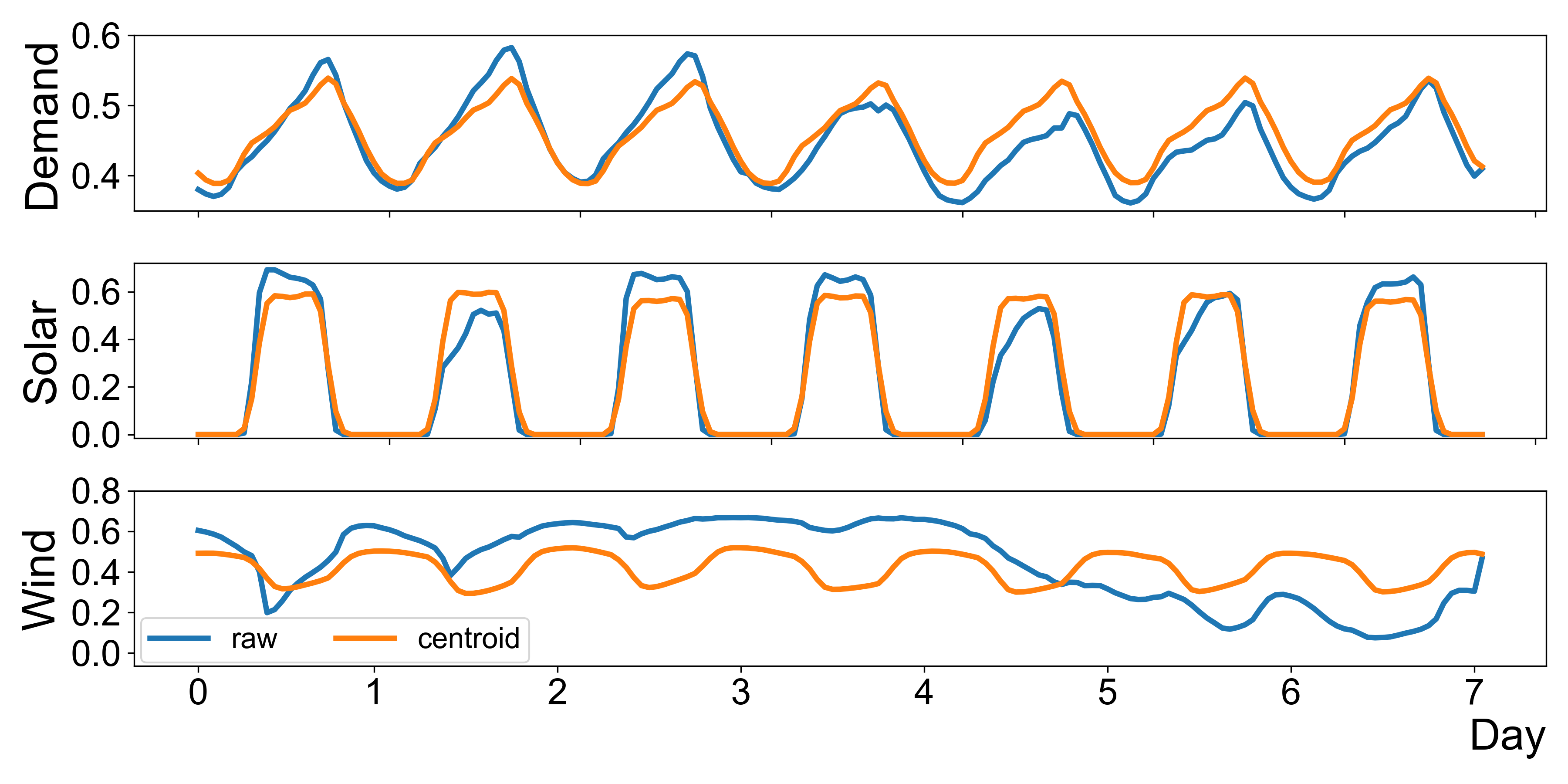}
    \caption{Representative week data of wind \& solar power and load curve}
    \label{centroid_smoothing}
\end{figure}

    After clustering, either the centroid points or the time slices closest to the centroids of each cluster can be used as the representative periods, with the number of slices in that cluster assigned as corresponding weights \cite{RN411}. Here we examine two settings: using the centroid (centroid option True) and using the time slice closest to the centroid (centroid option False). As shown in Figure \ref{centroid}, using the centroid increases the approximation error for all numbers of clusters. Mathematically, using the centroid means that the variation in the same cluster is averaged by other points. As shown in Figure \ref{centroid_smoothing}, the centroid point usually represents smoother profiles of wind, solar, and load, and those extreme values in the raw profiles will be lost if the centroid is used as the representative point. For example, the raw wind power data of the week closest to the centroid has a larger variation than the profile of the centroid, and the centroid profile gets smoother and more periodic after averaging. In ESOMs, the capacity expansion decisions are not only dependent on typical days/weeks but also extreme days/weeks \cite{young2014program}, and may be even more sensitive to the extreme scenarios which challenge the reliability of energy systems. Across representative periods, common and rare events are both appropriately selected and weighted through the clustering method, but within representative periods, averaged data can still smooth the variation especially in IRES, resulting in underestimation of required capacities like the wind power and storage capacities.


\section{Conclusions}
\label{Conclusions}
    In this article we proposed a model-adaptive clustering-based time aggregation method to better select representative periods for solving long-term ESOMs. New features are extracted based on ESOMs to improve approximation performance over traditional time aggregation. Two versions of an ESOM are designed to validate the proposed method, and the results show that the method can sufficiently reduce the approximation errors, compared to traditional time aggregation, while keeping the ESOM tractable. On the basis of clustering methods, deeper exploration into feature extraction and method settings are made to best exploit the performance of the proposed method. 
    
    By implementing this model-adaptive time aggregation method, the features from ESOMs are linked to each time period, allowing for a representation with lower dimensions. These features do help capture the characteristics of both variability in renewable sources and ESOMs. Only decision variables are used in the method validation, but there are lots of other features like dual variables of the ESOM. Exploring more relevant features from ESOMs is an interesting direction of future study.

\bibliography{Ref}{}

\begin{thebibliography}{10}
\providecommand{\url}[1]{#1}
\csname url@samestyle\endcsname
\providecommand{\newblock}{\relax}
\providecommand{\bibinfo}[2]{#2}
\providecommand{\BIBentrySTDinterwordspacing}{\spaceskip=0pt\relax}
\providecommand{\BIBentryALTinterwordstretchfactor}{4}
\providecommand{\BIBentryALTinterwordspacing}{\spaceskip=\fontdimen2\font plus
\BIBentryALTinterwordstretchfactor\fontdimen3\font minus
  \fontdimen4\font\relax}
\providecommand{\BIBforeignlanguage}[2]{{%
\expandafter\ifx\csname l@#1\endcsname\relax
\typeout{** WARNING: IEEEtran.bst: No hyphenation pattern has been}%
\typeout{** loaded for the language `#1'. Using the pattern for}%
\typeout{** the default language instead.}%
\else
\language=\csname l@#1\endcsname
\fi
#2}}
\providecommand{\BIBdecl}{\relax}
\BIBdecl

\bibitem{allen2019technical}
M.~Allen, P.~Antwi-Agyei, F.~Aragon-Durand, M.~Babiker, P.~Bertoldi, M.~Bind,
  S.~Brown, M.~Buckeridge, I.~Camilloni, A.~Cartwright \emph{et~al.},
  ``Technical summary: Global warming of 1.5° c. an ipcc special report on the
  impacts of global warming of 1.5° c above pre-industrial levels and related
  global greenhouse gas emission pathways, in the context of strengthening the
  global response to the threat of climate change, sustainable development, and
  efforts to eradicate poverty,'' 2019.

\bibitem{RN434}
M.~Da~Graça~Carvalho, ``Eu energy and climate change strategy,''
  \emph{Energy}, vol.~40, no.~1, pp. 19--22, 2012.

\bibitem{davidson2021policies}
M.~Davidson, V.~J. Karplus, D.~Zhang, and X.~Zhang, ``Policies and institutions
  to support carbon neutrality in china by 2060.'' \emph{Economics of Energy \&
  Environmental Policy}, vol.~10, no.~2, pp. 7--25, 2021.

\bibitem{sheet2021president}
F.~Sheet, ``President biden sets 2030 greenhouse gas pollution reduction target
  aimed at creating good-paying union jobs and securing us leadership on clean
  energy technologies,'' \emph{The White House}, 2021.

\bibitem{marcucci2017road}
A.~Marcucci, S.~Kypreos, and E.~Panos, ``The road to achieving the long-term
  paris targets: energy transition and the role of direct air capture,''
  \emph{Climatic Change}, vol. 144, no.~2, pp. 181--193, 2017.

\bibitem{wheeler2017carbon}
S.~M. Wheeler, ``A carbon-neutral california: Social ecology and prospects for
  2050 ghg reduction,'' \emph{Urban Planning}, vol.~2, no.~4, pp. 5--18, 2017.

\bibitem{RN366}
N.~L. Panwar, S.~C. Kaushik, and S.~Kothari, ``Role of renewable energy sources
  in environmental protection: A review,'' \emph{Renewable \& Sustainable
  Energy Reviews}, vol.~15, no.~3, pp. 1513--1524, 2011.

\bibitem{iea2020global}
U.~IEA, ``Global energy review 2020,'' \emph{Ukraine.[Online] https://www. iea.
  org/countries/ukraine [Accessed: 2020-09-10]}, 2020.

\bibitem{fridleifsson2001geothermal}
I.~B. Fridleifsson, ``Geothermal energy for the benefit of the people,''
  \emph{Renewable and sustainable energy reviews}, vol.~5, no.~3, pp. 299--312,
  2001.

\bibitem{bessa2014handling}
R.~Bessa, C.~Moreira, B.~Silva, and M.~Matos, ``Handling renewable energy
  variability and uncertainty in power systems operation,'' \emph{Wiley
  Interdisciplinary Reviews: Energy and Environment}, vol.~3, no.~2, pp.
  156--178, 2014.

\bibitem{golden2015curtailment}
R.~Golden and B.~Paulos, ``Curtailment of renewable energy in california and
  beyond,'' \emph{The Electricity Journal}, vol.~28, no.~6, pp. 36--50, 2015.

\bibitem{denholm2015overgeneration}
P.~Denholm, M.~O'Connell, G.~Brinkman, and J.~Jorgenson, ``Overgeneration from
  solar energy in california. a field guide to the duck chart,'' National
  Renewable Energy Lab.(NREL), Golden, CO (United States), Tech. Rep., 2015.

\bibitem{RN428}
S.~Koohi-Fayegh and M.~Rosen, ``A review of energy storage types, applications
  and recent developments,'' \emph{Journal of Energy Storage}, vol.~27, p.
  101047, 2020.

\bibitem{oudalov2007sizing}
A.~Oudalov, R.~Cherkaoui, and A.~Beguin, ``Sizing and optimal operation of
  battery energy storage system for peak shaving application,'' in \emph{2007
  IEEE Lausanne Power Tech}.\hskip 1em plus 0.5em minus 0.4em\relax IEEE, 2007,
  pp. 621--625.

\bibitem{stroe2016operation}
D.-I. Stroe, V.~Knap, M.~Swierczynski, A.-I. Stroe, and R.~Teodorescu,
  ``Operation of a grid-connected lithium-ion battery energy storage system for
  primary frequency regulation: A battery lifetime perspective,'' \emph{IEEE
  transactions on industry applications}, vol.~53, no.~1, pp. 430--438, 2016.

\bibitem{zhu2018optimization}
Y.~Zhu, C.~Liu, K.~Sun, D.~Shi, and Z.~Wang, ``Optimization of battery energy
  storage to improve power system oscillation damping,'' \emph{IEEE
  Transactions on Sustainable Energy}, vol.~10, no.~3, pp. 1015--1024, 2018.

\bibitem{ribeiro2001energy}
P.~F. Ribeiro, B.~K. Johnson, M.~L. Crow, A.~Arsoy, and Y.~Liu, ``Energy
  storage systems for advanced power applications,'' \emph{Proceedings of the
  IEEE}, vol.~89, no.~12, pp. 1744--1756, 2001.

\bibitem{ralon2017electricity}
P.~Ralon, M.~Taylor, A.~Ilas, H.~Diaz-Bone, and K.~Kairies, ``Electricity
  storage and renewables: Costs and markets to 2030,'' \emph{International
  Renewable Energy Agency: Abu Dhabi, UAE}, 2017.

\bibitem{alizadeh2011reliability}
B.~Alizadeh and S.~Jadid, ``Reliability constrained coordination of generation
  and transmission expansion planning in power systems using mixed integer
  programming,'' \emph{IET generation, transmission \& distribution}, vol.~5,
  no.~9, pp. 948--960, 2011.

\bibitem{haghighat2018bilevel}
H.~Haghighat and B.~Zeng, ``Bilevel mixed integer transmission expansion
  planning,'' \emph{IEEE Transactions on Power Systems}, vol.~33, no.~6, pp.
  7309--7312, 2018.

\bibitem{khodaei2012coordination}
A.~Khodaei, M.~Shahidehpour, L.~Wu, and Z.~Li, ``Coordination of short-term
  operation constraints in multi-area expansion planning,'' \emph{IEEE
  Transactions on Power Systems}, vol.~27, no.~4, pp. 2242--2250, 2012.

\bibitem{castro2018expanding}
P.~M. Castro, I.~E. Grossmann, and Q.~Zhang, ``Expanding scope and
  computational challenges in process scheduling,'' \emph{Computers \& Chemical
  Engineering}, vol. 114, pp. 14--42, 2018.

\bibitem{RN429}
V.~Oree, S.~Z. Sayed~Hassen, and P.~J. Fleming, ``Generation expansion planning
  optimisation with renewable energy integration: A review,'' \emph{Renewable
  and Sustainable Energy Reviews}, vol.~69, pp. 790--803, 2017.

\bibitem{RN185}
H.~Teichgraeber, C.~P. Lindenmeyer, N.~Baumgartner, L.~Kotzur, D.~Stolten,
  M.~Robinius, A.~Bardow, and A.~R. Brandt, ``Extreme events in time series
  aggregation: A case study for optimal residential energy supply systems,''
  \emph{Applied Energy}, vol. 275, p. 115223, 2020.

\bibitem{RN192}
X.~Zhang and D.~J. Hil, ``Hierarchical temporal and spatial clustering of
  uncertain and time-varying load models,'' 2020.

\bibitem{haydt2011relevance}
G.~Haydt, V.~Leal, A.~Pina, and C.~A. Silva, ``The relevance of the energy
  resource dynamics in the mid/long-term energy planning models,''
  \emph{Renewable energy}, vol.~36, no.~11, pp. 3068--3074, 2011.

\bibitem{pina2013high}
A.~Pina, C.~A. Silva, and P.~Ferr{\~a}o, ``High-resolution modeling framework
  for planning electricity systems with high penetration of renewables,''
  \emph{Applied Energy}, vol. 112, pp. 215--223, 2013.

\bibitem{welsch2015supporting}
M.~Welsch, M.~Howells, M.~R. Hesamzadeh, B.~{\'O}~Gallach{\'o}ir, P.~Deane,
  N.~Strachan, M.~Bazilian, D.~M. Kammen, L.~Jones, G.~Strbac \emph{et~al.},
  ``Supporting security and adequacy in future energy systems: The need to
  enhance long-term energy system models to better treat issues related to
  variability,'' \emph{International Journal of Energy Research}, vol.~39,
  no.~3, pp. 377--396, 2015.

\bibitem{nahmmacher2016carpe}
P.~Nahmmacher, E.~Schmid, L.~Hirth, and B.~Knopf, ``Carpe diem: A novel
  approach to select representative days for long-term power system modeling,''
  \emph{Energy}, vol. 112, pp. 430--442, 2016.

\bibitem{RN330}
M.~Zatti, M.~Gabba, M.~Freschini, M.~Rossi, A.~Gambarotta, M.~Morini, and
  E.~Martelli, ``k-milp: A novel clustering approach to select typical and
  extreme days for multi-energy systems design optimization,'' \emph{Energy},
  vol. 181, pp. 1051--1063, 2019.

\bibitem{RN420}
R.~Domínguez and S.~Vitali, ``Multi-chronological hierarchical clustering to
  solve capacity expansion problems with renewable sources,'' \emph{Energy},
  vol. 227, p. 120491, 2021.

\bibitem{poncelet2020unit}
K.~Poncelet, E.~Delarue, and W.~D’haeseleer, ``Unit commitment constraints in
  long-term planning models: Relevance, pitfalls and the role of assumptions on
  flexibility,'' \emph{Applied Energy}, vol. 258, p. 113843, 2020.

\bibitem{RN293}
W.~W. Tso, C.~D. Demirhan, C.~F. Heuberger, J.~B. Powell, and E.~N.
  Pistikopoulos, ``A hierarchical clustering decomposition algorithm for
  optimizing renewable power systems with storage,'' \emph{Applied Energy},
  vol. 270, p. 115190, 2020.

\bibitem{RN411}
L.~Reichenberg, A.~S. Siddiqui, and S.~Wogrin, ``Policy implications of
  downscaling the time dimension in power system planning models to represent
  variability in renewable output,'' \emph{Energy}, vol. 159, pp. 870--877,
  2018.

\bibitem{Zhang_Adaptive-clustering-for-ESOM_2021}
\BIBentryALTinterwordspacing
Y.~Zhang and G.~He, ``{Adaptive-clustering-for-ESOM},'' 2021. [Online].
  Available: \url{https://github.com/Betristor/Adaptive-clustering-for-ESOM}
\BIBentrySTDinterwordspacing

\bibitem{RN509}
\BIBentryALTinterwordspacing
A.~Garcia-Cerezo, L.~Baringo, and R.~Garcia-Bertrand, ``Representative days for
  expansion decisions in power systems,'' \emph{Energies}, vol.~13, no.~2,
  2020. [Online]. Available: \url{https://www.mdpi.com/1996-1073/13/2/335}
\BIBentrySTDinterwordspacing

\bibitem{RN195}
M.~Gamst, S.~Buchholz, and D.~Pisinger, ``Time aggregation techniques applied
  to a capacity expansion model for real-life sector coupled energy systems,''
  Conference Proceedings.

\bibitem{RN419}
N.~González-Cabrera, J.~Ortiz-Bejar, A.~Zamora-Mendez, and M.~R.
  Arrieta~Paternina, ``On the improvement of representative demand curves via a
  hierarchical agglomerative clustering for power transmission network
  investment,'' \emph{Energy}, vol. 222, p. 119989, 2021.

\bibitem{RN424}
T.~M. Kodinariya and P.~R. Makwana, ``Review on determining number of cluster
  in k-means clustering,'' \emph{International Journal}, vol.~1, no.~6, pp.
  90--95, 2013.

\bibitem{young2014program}
D.~Young \emph{et~al.}, ``Program on technology innovation: Us-regen model
  documentation 2014,'' \emph{EPRI, CA}, 2014.

\end{thebibliography}
\bibliographystyle{IEEEtran}
\end{document}